\documentclass{article}

\usepackage[T1]{fontenc}
\usepackage{amsmath}
\usepackage{amssymb}
\usepackage{amsthm}
\usepackage{mathrsfs}
\usepackage{pxfonts}
\usepackage{sectsty}
\usepackage{enumitem}
\usepackage{microtype}
\usepackage{xcolor}
\usepackage[hidelinks]{hyperref}
\usepackage{url}
\hypersetup{
 pdftitle={On some stochastic evolutions on Laguerre-Polya functions},
 pdfauthor={Theodoros Assiotis and Fengyi Li}
}

\allowdisplaybreaks
\sectionfont{\scshape\centering\fontsize{11}{14}\selectfont}
\subsectionfont{\scshape\fontsize{11}{14}\selectfont}

\usepackage{fancyhdr}
\newcommand\shorttitle{On Laguerre--P\'olya stochastic evolutions}
\newcommand\authors{T. Assiotis and F. Li}
\newtheorem{thm}{Theorem}[section]
\newtheorem{lem}[thm]{Lemma}
\newtheorem{prop}[thm]{Proposition}
\newtheorem{cor}[thm]{Corollary}
\newtheorem{defn}[thm]{Definition}
\newtheorem{rmk}[thm]{Remark}
\theoremstyle{plain}

\newcommand{\E}{\mathbb E}
\newcommand{\R}{\mathbb R}
\newcommand{\N}{\mathbb N}
\newcommand{\dd}{\mathrm{d}}
\newcommand{\defeq}{\overset{\mathrm{def}}{=}}
\newcommand{\1}{\mathbf 1}
\newcommand{\Law}{\mathsf{Law}}
\newcommand{\Tr}{\operatorname{Tr}}
\newcommand{\C}{\mathbb C}

\title{\large\bf ON SOME STOCHASTIC EVOLUTIONS ON LAGUERRE--P\'OLYA FUNCTIONS}
\author{\small THEODOROS ASSIOTIS AND FENGYI LI}
\date{}

\newcommand{\resultspace}[1][5]{\par\begingroup
 \dimen0=\pagegoal\advance\dimen0 by -\pagetotal
 \ifdim\dimen0<#1\baselineskip\newpage\fi\endgroup}

\begin{document}

\maketitle

\begin{abstract}
We prove uniqueness in law for three stochastic evolutions on normalised
Laguerre--P\'olya functions arising as scaling limits of matrix
characteristic-polynomial dynamics.  In particular, this implies that Ahn's universal
Brownian Gibbs line ensemble arising from random matrix products has the law of the logarithms of the
reciprocal zeros of the solution, unique in law, to a simple stochastic
equation on entire functions started from \(\mathrm e^{-z}\).
 The method of proof relies on some integrable structure of the problem in the special form of a new duality with finite-dimensional diffusions which  is closely related to intertwinings.
\end{abstract}

\section{Introduction}
\label{sec:introduction}

\subsection{Setting and motivation}

Certain infinite-dimensional particle dynamics arise naturally as microscopic
limits of eigenvalue and singular-value evolutions of random matrices.
Their description by infinite systems of stochastic differential equations
(ISDEs) with logarithmic interaction raises basic questions of existence and
uniqueness given the choice of initial state.  Important approaches to these dynamics include
the determinantal constructions of Katori--Tanemura
\cite{KatoriTanemura}, the Dirichlet-form and ISDE theory of Osada and
Osada--Tanemura \cite{OsadaRandomMatrices,OsadaTanemura}, and Tsai's construction of
strong solutions to the infinite Dyson system \cite{Tsai}; see also the authors' paper \cite{assiotis2026infinite}.

In recent  work \cite{AssiotisMirsajjadiISDE},
the infinite-particle limit of the so-called Rider--Valk\'o dynamics from \cite{RiderValkoMatrixDufresne} was constructed and shown that
its logarithmic-interaction ISDE,
\begin{equation}
\label{eq:Rider-Valko-ISDE}
 \dd\mathsf x_i(t)
 =\mathsf x_i(t)\,\dd\mathsf B_i(t)
 -\frac{\nu}{2}\mathsf x_i(t)\,\dd t
 +\sum_{\substack{j=1\\j\neq i}}^{\infty}
   \frac{\mathsf x_i(t)\mathsf x_j(t)}
        {\mathsf x_i(t)-\mathsf x_j(t)}\,\dd t,
 \qquad i\in\N,
\end{equation}
with $(\mathsf B_i)_{i\geq1}$ independent standard Brownian
motions and $\mathsf x_1(t)>\mathsf x_2(t)>\cdots>0$, admits distinct (weak) solution laws from the
same initial particle configuration. This provided
the first concrete examples of such non-uniqueness for an ISDE with
logarithmic interaction.  To retain  additional information that one may hope to need for
a characterisation, the full evolution of the
limiting reverse characteristic polynomial was subsequently considered in \cite{AssiotisMirsajjadi}, also for two other models: the limiting singular values of multiplicative Brownian motion on the general linear group \(\mathscr{GL}_N(\C)\) \cite{NorrisRogersWilliams,JonesOConnell,Ahn} and the dynamical Cauchy or Hua-Pickrell model \cite{AssiotisHP,AssiotisBougerol,AristaDemni}.
The resulting state spaces consist of (subsets of) normalised Laguerre--P\'olya
functions, which arise as locally uniform limits of polynomials with real
zeros and have the canonical form
\[
 \begin{gathered}
 \mathrm{E}\left(z\right)=\mathrm e^{-\gamma z-qz^2/2}
 \prod_{j=1}^{\infty}\mathrm e^{y_jz}\left(1-y_jz\right),\\
 \gamma\in\R,\quad q\geq0,\quad y_j\in\R,\quad
 \sum_{j=1}^{\infty}y_j^2<\infty.
 \end{gathered}
\]
The exponential factors record information that is invisible to the zero
set alone.  This enlarged description distinguishes evolutions whose
particle coordinates may solve the same ISDE.

Here, as predicted in \cite{AssiotisMirsajjadi}, we prove that the three entire-function valued equations of
\cite{AssiotisMirsajjadi} do characterise their solutions in law.
Such an intrinsic characterisation gives a way to identify subsequential
limits and is therefore useful for convergence and universality problems.
As a corollary, we obtain a stochastic analytic characterisation of Ahn's
Brownian Gibbs line ensemble \cite{Ahn}: its curves are the logarithms of the
reciprocal zeros of a solution started from \(\mathrm e^{-z}\).
This ensemble is expected to describe a broad universality class for
products of many large random matrices: in particular Ahn \cite{Ahn} proves convergence
of the leading curves in finite-dimensional distributions for a class of
products of independent right-unitarily invariant matrices, after
centering and a time change, when the number of factors is comparable
to the matrix size.  

The proof here relies in an essential way on the integrable structure of the problem in terms of a new duality with finite-dimensional diffusions. Based on this duality the arguments are not very complicated and rather elegant. Nevertheless, it would be desirable to have a more robust/less integrable approach to problems of this kind. Finally, a different way to restore uniqueness to the ISDEs has been suggested in \cite{AssiotisMirsajjadiISDE,AssiotisMirsajjadi} by considering so-called "rigid-path" solutions, see \cite{AssiotisMirsajjadiISDE,AssiotisMirsajjadi} for a more detailed discussion. Some preliminary progress with this viewpoint in mind is made in \cite{MirsajjadiDeterminantal}. It would be interesting to understand this other type of characterisation.

\subsection{Main results}

We now move to state our main results precisely. Two of the stochastic evolutions live in the following positive class of Laguerre-P\'olya functions
\begin{equation}
 \label{eq:intro-positive-state}
 \mathscr E_+=\left\{\mathrm e^{-dz}\prod_{j=1}^{\infty}\left(1-x_jz\right):
 d\geq0,\ x_1\geq x_2\geq\cdots\geq0,\
 \sum_{j=1}^{\infty}x_j<\infty\right\}.
\end{equation}
They are the edge-limit equation associated with Brownian motion on
\(\mathscr{GL}_N(\C)\) in
\cite{AssiotisMirsajjadi} and the limiting dynamics of the Rider--Valk\'o equation from the matrix
Dufresne identity framework \cite{RiderValkoMatrixDufresne}.  The third evolves on
the full normalised two-sided class of Laguerre-P\'olya functions
\begin{equation}
 \label{eq:intro-HP-state}
 \mathscr E
 =\left\{\mathrm e^{-\gamma z-qz^2/2}
   \prod_{j=1}^{\infty}\mathrm e^{y_jz}\left(1-y_jz\right):
   \gamma\in\R,\ q\geq0,\ \boldsymbol y=\left(y_j\right)\in\ell^2\left(\N;\R\right)\right\},
\end{equation}
where \(\ell^2(\N;\R)\) is the space of square-summable real sequences,
and the nonzero \(y_j\)'s are put in the signed canonical order fixed
below.  Both spaces carry the topology of locally uniform convergence.
For \(g\in\mathscr E\), put
\(\gamma(g)=-g'(0)=\gamma\) and
\(\delta(g)=q+\sum_{j=1}^{\infty}y_j^2\); for \(g\in\mathscr E_+\), put
\(\gamma(g)=-g'(0)=d+\sum_{j=1}^{\infty}x_j\).

For every normalised entire function \(g\), define the analytic "covariance"
kernel
\begin{equation}
 \label{eq:intro-Q}
 \mathcal Q_g\left(z,w\right)
 =\frac{zw}{z-w}
 \left(g\left(z\right)g'\left(w\right)-g\left(w\right)g'\left(z\right)\right),\qquad
 \mathcal Q_g^{\mathrm{HP}}=2\mathcal Q_g,
\end{equation}
with the diagonal values given by analytic continuation.  The positive
equations have drift
\begin{equation}
 \label{eq:intro-drift}
 \left(\mathcal L_{\theta,\kappa}g\right)\left(z\right)
 =\frac{\theta}{2}zg'\left(z\right)-\frac12z^2g''\left(z\right)-\frac{\kappa}{2}zg\left(z\right),
 \qquad\theta\in\R,\quad\kappa\in\{0,1\}.
\end{equation}
The \(\mathscr{GL}_N(\C)\) edge-limit equation has \(\kappa=0\); the
Rider--Valk\'o equation with parameter \(\nu\) has
\((\theta,\kappa)=(-\nu,1)\).  For the Hua--Pickrell parameter
\(\mathfrak s=a+\mathrm i b\), the drift is given by,
\begin{equation}
 \label{eq:intro-HP-drift}
 \left(\mathcal L_{a,b}^{\mathrm{HP}}g\right)\left(z\right)
 =-\left(2bz+z^2\right)g\left(z\right)-2azg'\left(z\right)-z^2g''\left(z\right).
\end{equation}

\begin{defn}
 Set \(\left(\mathcal L^+,\mathcal Q^+\right)
=(\mathcal L_{\theta,\kappa},\mathcal Q)\) and
\(\mathcal L^{\mathrm{HP}}=\mathcal L_{a,b}^{\mathrm{HP}}\).
For \(\diamond\in\{+,\mathrm{HP}\}\), a weak solution
\(\boldsymbol{\Phi}\) is defined on a filtered probability space
\[
(\Omega,\mathscr F,(\mathscr F_t)_{t\geq0},\mathbb P)
\]
satisfying the usual conditions, and satisfies
\begin{equation}
 \label{eq:intro-martingale}
 M_t^{\boldsymbol{\Phi}}\left(z\right)
 =\boldsymbol{\Phi}_t\left(z\right)
  -\boldsymbol{\Phi}_0\left(z\right)
  -\int_0^t
    \left(\mathcal L^\diamond\boldsymbol{\Phi}_r\right)
    \left(z\right)\,\dd r,
 \qquad z\in\C,
\end{equation}
where the processes \(M^{\boldsymbol{\Phi}}(z)\) are continuous
local martingales on one filtration and have complex-bilinear brackets
\begin{equation}
 \label{eq:intro-bracket}
 \left[M^{\boldsymbol{\Phi}}\left(z\right),
       M^{\boldsymbol{\Phi}}\left(w\right)\right]_t
 =\int_0^t
   \mathcal Q_{\boldsymbol{\Phi}_r}^\diamond
   \left(z,w\right)\,\dd r.
\end{equation}
\end{defn}

The above martingale formulation can be written in equivalent form as an evolution equation with respect to a cylindrical Brownian motion. More precisely, on a possibly enlarged probability space, \eqref{eq:intro-martingale}, \eqref{eq:intro-bracket} it can be written as
\[
\dd\boldsymbol{\Phi}_t=\mathcal{L}_\diamond\boldsymbol{\Phi}_t\dd t+\Sigma_\diamond(\boldsymbol{\Phi}_t)\,\dd\mathsf{W}_t,
\]
where \(\mathsf{W}\) is a real cylindrical Brownian motion on
\(\ell^2(\mathbb N_0;\mathbb R)\), with independent standard real
coordinate Brownian motions \((\mathsf B_j)_{j\ge0}\). Writing
\begin{equation*}
    \begin{aligned}
    \boldsymbol{\Phi}_t (z)&= \mathrm{e}^{-\mathsf d_t z}\prod_{j=1}^\infty(1-\mathsf x_j(t)z),\quad\text{if } \boldsymbol{\Phi}_t\in \mathscr E_+,\\
    \boldsymbol{\Phi}_t (z)&=\mathrm{e}^{-\boldsymbol\gamma_t z-\mathsf q_tz^2/2}\prod_{j=1}^\infty\mathrm{e}^{\mathsf y_j(t)z}(1-\mathsf y_j(t)z)\quad\text{if } \boldsymbol{\Phi}_t\in \mathscr E,
    \end{aligned}
\end{equation*}
respectively, the noise coefficient \(\Sigma_\diamond\) is specified through evaluations by
\[
\begin{aligned}
\bigl(\Sigma_+(\boldsymbol{\Phi}_t)\,\dd\mathsf W_t\bigr)(z)
&=-\sum_{j\ge1}
  \frac{z \mathsf x_j(t)\boldsymbol{\Phi}_t(z)}{1-\mathsf x_j(t)z}\,\dd\mathsf B_j(t),\\
\bigl(\Sigma_{\mathrm{HP}}(\boldsymbol{\Phi}_t)\,\dd\mathsf W_t\bigr)(z)
&=-\sqrt{2\mathsf q_t}\,z\boldsymbol{\Phi}_t(z)\,\dd\mathsf B_0(t)-\sqrt{2}\sum_{j\ge1}
  \frac{z \mathsf y_j(t)\boldsymbol{\Phi}_t(z)}{1-\mathsf y_j(t)z}\,\dd\mathsf B_j(t).
\end{aligned}
\]
Note that, the \(\mathsf x_j(t)\) and \(\mathsf y_j(t)\) are the canonical positive and signed reciprocal-root coordinates of \(\boldsymbol{\Phi}_t\), respectively, and \(\boldsymbol \gamma_t=\gamma(\boldsymbol{\Phi}_t)\) and
\(\mathsf q_t=q(\boldsymbol{\Phi}_t)\). The apparent poles are removable, and the series are interpreted as localised stochastic
integrals. For details on cylindrical stochastic integration and Brownian representation, see \cite[Section~4.4]{HairerSPDEs} and \cite[Theorem~2]{Ondrejat2005}. This formulation is closer to the stochastic partial differential equations literature and would be interesting to develop further, but we will not need this here.

The so-called boundary processes of
\cite{AssiotisHP,AssiotisMirsajjadiISDE} are constructed by the method of intertwiners, which builds a Markov
semigroup on a projective boundary from finite-dimensional semigroups
consistent under Markov kernels; see
\cite{BorodinOlshanskiGTBoundary,BorodinOlshanskiGTProcesses,
BorodinOlshanskiThoma}.
Their entire-function limits in
\cite[Theorems~1.4 and~1.6]{AssiotisMirsajjadi}
provide weak solutions $\boldsymbol{\Phi}_t$ to the martingale problems we defined.
Our main theorem then identifies their laws uniquely from the closed
martingale problems above.

\begin{thm}
\label{thm:intro-uniqueness}
For every \(\theta\in\R\), \(\kappa\in\{0,1\}\), and deterministic
\(g\in\mathscr E_+\), the martingale problem
\eqref{eq:intro-martingale}--\eqref{eq:intro-bracket} with
\(\mathcal L^\diamond=\mathcal L_{\theta,\kappa}\) and
\(\mathcal Q^\diamond=\mathcal Q\) has a unique in law weak solution in \(C(\R_+;\mathscr E_+)\).  The analogous problem has a unique in law weak solution on \(C(\R_+;\mathscr E)\), for
every \(a,b\in\R\) and deterministic \(g\in\mathscr E\),
with \(\mathcal L^\diamond=\mathcal L_{a,b}^{\mathrm{HP}}\) and
\(\mathcal Q^\diamond=\mathcal Q^{\mathrm{HP}}\).
\end{thm}

The uniqueness theorem above also characterises Ahn's \(\mathscr{GL}_N(\C)\) Brownian motion singular-value limit through
the zeros of one stochastic entire function.  We first fix its scaling
conventions.

\begin{defn}
\label{def:ahn-line-ensemble}
Let \(\mathbf Y_N\) be Brownian motion on \(\mathscr{GL}_N(\C)\) given by
the It\^o equation
\[
 \dd\mathbf Y_N\left(t\right)=\mathbf Y_N\left(t\right)\dd\mathbf B_N\left(t\right),\qquad
 \mathbf Y_N\left(0\right)=\mathbf I_N,
\]
where \(\mathbf B_N=\mathbf B_N^{(1)}+\mathrm i\mathbf B_N^{(2)}\) and
all entries of the real matrices \(\mathbf B_N^{(1)},\mathbf B_N^{(2)}\)
are independent standard Brownian motions.  If
\(\sigma_1^{(N)}(t)\geq\cdots\geq\sigma_N^{(N)}(t)>0\) are its singular
values, put
\begin{equation}
 \xi_i^{\left(N\right)}\left(t\right)=2\log\sigma_i^{\left(N\right)}\left(\frac t4\right)
 -\frac{Nt}{2}-\log N,\qquad 1\leq i\leq N.
 \label{eq:ahn-scaling}
\end{equation}
Ahn's line ensemble \(\boldsymbol{\xi}=(\xi_i(t))_{i\in\N,\,t>0}\)
is their limit: for every \(k\geq1\) and \(0<\varepsilon<T\), the first
\(k\) curves converge in distribution in \(C([\varepsilon,T];\R^k)\)
\cite[Theorem~1.3]{Ahn}.
\end{defn}

Our general weak uniqueness result identifies this ensemble from the initial state
entire function \(\mathrm e^{-z}\).

\begin{cor}
\label{cor:ahn-spde-identification}
Let \(\boldsymbol{\Phi}\) solve the positive martingale problem of
Theorem~\ref{thm:intro-uniqueness} with \((\theta,\kappa)=(0,0)\) and
\(\boldsymbol{\Phi}_0(z)=\mathrm e^{-z}\).  Its ordered positive zeros
\(r_1(t)\leq r_2(t)\leq\cdots\), counted with multiplicity, form an
infinite continuous family for \(t>0\), and
\begin{equation}
 \Law\left(\left(-\log r_i\left(t\right)\right)_{i\in\N,\,t>0}\right)
 =\Law\left(\boldsymbol{\xi}\right).
 \label{eq:ahn-zero-identification}
\end{equation}
\end{cor}

\begin{proof}
Set \(x_i^{(N)}(t)=\mathrm e^{\xi_i^{(N)}(t)}\).  The reverse polynomials
\(\prod_{i=1}^{N}(1-zx_i^{(N)}(t))\) start from
\((1-z/N)^N\to \mathrm e^{-z}\). By \cite[Theorem~1.4 and Eq.~(12)]{AssiotisMirsajjadi}, these processes converge jointly in path space with the scaled squared singular values. The limiting entire-function-valued process solves the present martingale problem, and its canonical reciprocal roots coincide with the limiting particle coordinates. By \eqref{eq:ahn-scaling} and
Ahn's convergence, those coordinate paths have law
\((\mathrm e^{\xi_i(t)})_{i\in\N,\,t>0}\).  Theorem~\ref{thm:intro-uniqueness}
identifies the entire-function path law with that of \(\boldsymbol{\Phi}\),
proving the assertion.
\end{proof}

\subsection{On the proof}

The whole proof rests on a single idea which exploits the integrable structure of the problem. Our use of integrability is through finite-dimensional duality.  The
reciprocal products
\(H_n^+(g,\mathbf s)=\prod_{j=1}^{n} g(-s_j)^{-1}\) and
\(H_n^{\mathrm{HP}}(g,\mathbf u)=\prod_{j=1}^{n} g(\mathrm i u_j)^{-1}\)
are bounded in modulus by one and determine probability laws on the
respective spaces of entire functions.  Applying It\^o's formula to
these probes turns the nonlinear covariance into a finite-dimensional
diffusion generator and a multiplication operator.  For each \(n\in \mathbb{N}\) the dual diffusions \(\boldsymbol{\mathsf S}=(\mathsf S_i)_{i=1}^n\)
have the common form
\[
 \dd\mathsf S_i=\sqrt c\,\mathsf S_i\,\dd\mathsf B_i+
 \left(\beta\mathsf S_i+c\sum_{\substack{j=1\\j\neq i}}^{n}
 \frac{\mathsf S_i\mathsf S_j}{\mathsf S_i-\mathsf S_j}\right)\dd t,
 \qquad 1\leq i\leq n,
\]
where the Brownian motions are independent and
\((c,\beta)=(1,\theta/2)\) in the positive case and
\((c,\beta)=(2,-2a)\) in the Hua--Pickrell case.  The first diffusion
lives in \(0<s_1<\cdots<s_n\); the second lives in a signed chamber
\(u_1<\cdots<u_\ell<0<u_{\ell+1}<\cdots<u_n\).
Both are conservative, avoid collisions, and preserve signs.

This essentially gives the main step of the proof which is the conditional duality of
Proposition~\ref{prop:conditional-duality}: for \(\diamond=+\) or
\(\mathrm{HP}\),
\[
 \E\left[H_n^\diamond\left(\boldsymbol{\Phi}_t,\mathbf v\right)
             \mid\mathscr F_r\right]
 =\E_{\mathbf v}^\diamond\left[
 \mathrm e^{\int_0^{t-r}C_\diamond\left(\boldsymbol{\mathsf S}_q\right)\,\dd q}
 H_n^\diamond\left(\boldsymbol{\Phi}_r,\boldsymbol{\mathsf S}_{t-r}\right)\right],
 \qquad 0\leq r\leq t,
\]
where \(C_+(\mathbf s)=-\frac\kappa2\sum_{i=1}^{n} s_i\),
\(C_{\mathrm{HP}}(\mathbf u)=2\mathrm i b\sum_{i=1}^{n} u_i-\sum_{i=1}^{n} u_i^2\),
and the expectation on the right is over an independent dual diffusion \(\boldsymbol{\mathsf S}\).
 Then, localisation, together
with conditioning on bounded initial-coordinate sets in the
Hua--Pickrell case, proves this identity for arbitrary initial laws
without moment assumptions.  Since the probes determine laws,
iteration of the conditional identity gives uniqueness of the entire
path law and hence the Markov semigroups.

\subsection{Quantitative convergence to equilibrium}
The duality methods developed here also provide a starting point for
studying quantitative convergence to equilibrium towards stochastic
zeta functions; see
\cite{ChhaibiNajnudelNikeghbali,ValkoViragZeta,
LiValkoCircularJacobi,AssiotisEntireFunctions}
for related constructions and limits of characteristic polynomials.
In particular, it is natural to ask whether convergence from fixed
initial states is exponentially fast in Wasserstein distances induced
by bounded metrics compatible with the locally uniform topology.
Combining the duality viewpoint with coupling arguments should provide
a route to such estimates. In the Rider--Valk\'o and real-parameter
Hua--Pickrell cases, the equilibrium limits have the form
\[
 \boldsymbol{\Phi}_t
 \xrightarrow[t\to\infty]{\mathrm{law}}
 \boldsymbol{\Phi}_\infty,
 \qquad
 \boldsymbol{\Phi}_\infty(z)=
 \begin{cases}
 \displaystyle
 \prod_{\xi\in\mathsf{Bes}_\nu}
       \left(1-\frac{4z}{\xi}\right),
 & \text{Rider--Valk\'o},\quad \nu>-1,
 \\[1.2em]
 \displaystyle
 \operatorname{pv}\prod_{\xi\in\mathcal Z_a}
       \left(1-\frac{z}{\xi}\right),
 & \text{Hua--Pickrell},\quad a>-\frac12,\ b=0,
 \end{cases}
\]
with convergence in the locally uniform topology
\cite{AssiotisMirsajjadiISDE,AssiotisMirsajjadi}.
Here \(\mathsf{Bes}_\nu\) is the Bessel determinantal point process
in the Tracy--Widom normalization, with hard-edge variable
\(4Nx\) for the Laguerre weight \(x^\nu\mathrm e^{-x}\)
\cite{TracyWidomBessel,ForresterLogGases}, and \(\mathcal Z_a\)
is the reciprocal image of the real-parameter Hua--Pickrell boundary
point process \cite{BorodinOlshanski,Qiu}.
The principal value is taken through symmetric truncations
\(|\xi|<N^2\). At \(a=b=0\), the zeros form the sine determinantal
point process, up to a deterministic spatial rescaling
\cite{AssiotisHP}.
Quantitative relaxation to this equilibrium is particularly interesting
in view of Suzuki's result \cite{SuzukiNoSpectralGap}: the unlabelled
infinite Dyson Brownian motion at \(\beta=2\), which also has sine
equilibrium, has no spectral gap and hence admits no uniform
exponential convergence estimate in \(\mathrm L^2\) of its invariant measure.

\paragraph{Acknowledgements} Fengyi Li is grateful to the School of Mathematics of the University of Edinburgh for financial support through a PhD studentship.

\paragraph{AI use statement} Motivated by the use of duality to prove uniqueness for martingale formulations of stochastic partial differential equations, and the fact that dualities are known to be intimately related to intertwinings, we used AI tools (ChatGPT 5.6 Sol) to search for and find the precise form of the duality in Proposition \ref{prop:conditional-duality} a few months ago. Once this duality is found the rest of the argument is routine. We were also using AI tools (ChatGPT 5.5, 5.6 Sol) prior to this to investigate other approaches (which led to results) but the arguments were much less elegant to the duality route taken here. We used ChatGPT 6 to assist the writeup with rough drafts following instructions and outlines from the authors. We then revised the proofs and prepared the final manuscript and we are solely responsible for its contents.

\section{State spaces and martingale problems}
\label{sec:state-spaces}

\subsection{State-space preliminaries}

Throughout, we use the state spaces, drift operators and covariance
kernels introduced in Section~\ref{sec:introduction}.  A weak solution
with initial law \(\mu\) is an adapted process with continuous paths in
the corresponding state space, defined on a filtered probability space
satisfying the usual conditions, with initial distribution \(\mu\), such
that the evaluation processes in \eqref{eq:intro-martingale} are
continuous local martingales with the covariations
\eqref{eq:intro-bracket}, all relative to the same filtration.
Covariations are extended complex-bilinearly.  A solution from \(g\)
means one with initial law \(\delta_g\).

For \(g\in\mathscr E_+\), write \(\mathbf x=(x_j)_{j\geq1}\) for its
ordered reciprocal-root coordinates.  For
\(g\in\mathscr E\), let \(\boldsymbol\alpha^+\) and
\(\boldsymbol\alpha^-\) be the decreasing lists of positive and negative
reciprocal-root magnitudes, padded with zeros, and fix the signed enumeration
\begin{equation}
 \label{hp:eq:signed-root-enumeration}
 y_{2j-1}=\alpha_j^+,\qquad y_{2j}=-\alpha_j^-,\qquad j\geq1.
\end{equation}
The coordinates satisfy
\begin{equation}
 \label{hp:eq:coordinates-derivatives}
 \gamma(g)=-g'(0),\qquad
 \delta(g)=\gamma(g)^2-g''(0)\quad(g\in\mathscr E).
\end{equation}

The canonical products identify \(\mathscr E_+\) and
\(\mathscr E\) homeomorphically with their boundary spaces
in the coordinates \((\mathbf x,\gamma)\) and
\((\boldsymbol\alpha^+,\boldsymbol\alpha^-,\gamma,\delta)\), respectively,
equipped with coordinatewise convergence; see
\cite[Remarks~2.6--2.7 and the proof of Proposition~5.6]{AssiotisEntireFunctions}.
Both spaces are Polish, and the ordered reciprocal-root coordinates are
continuous.  The "dust" coordinates
\(d=\gamma-\sum_{j=1}^{\infty}x_j\) and \(q=\delta-\sum_{j=1}^{\infty}y_j^2\) are Borel measurable.
For every \(R<\infty\), the sets
\begin{equation}
 \label{hp:eq:compact-sublevels}
 \begin{split}
 K_R&=\{g\in\mathscr E_+:\gamma(g)\leq R\},\\
 \mathscr K_R&=\{g\in\mathscr E:
                         \gamma(g)^2+\delta(g)\leq R\}
 \end{split}
\end{equation}
are compact, since the corresponding coordinate sublevels are compact.

Continuous weak solutions from every deterministic state are supplied by
\cite[Theorems~1.4 and~1.6]{AssiotisMirsajjadi}.  Theorem~1.4 covers
\((\theta,\kappa)=(\theta,0)\); Theorem~1.6 covers
\((\theta,\kappa)=(-\nu,1)\) and Hua--Pickrell with
\(\mathfrak s=a+\mathrm i b\), and taking \(\nu=-\theta\) gives the
remaining positive cases.  Their initial-state hypotheses follow from
the density of normalised polynomials with positive, respectively real,
zeros and the equivalence of locally uniform convergence and convergence
of the boundary coordinates; see
\cite[Proposition~2.5 and Remark~2.7]{AssiotisEntireFunctions}.

\subsection{Hua--Pickrell covariance and bounded probes}

We first introduce reciprocal evaluations on the imaginary axis, which provide bounded observables for the duality and tightness arguments.
\begin{defn}
\label{hp:def:reciprocal-probes}
For $g\in\mathscr E$, define
\begin{equation}
 \label{hp:eq:h-def}
 h_g^{\mathrm{HP}}\left(u\right)\defeq\frac1{g\left(\mathrm i u\right)},\qquad u\in\R.
\end{equation}
For fixed \(u\), this is continuous in \(g\) for locally uniform
convergence.  If \(u\ne0\), the denominator is nonzero because every zero
of \(g\) is real; the case \(u=0\) is immediate.
\end{defn}

The characteristic-function representation supplies the derivative and
modulus bounds needed to justify duality and control state-space tails.
\begin{lem}
\label{hp:lem:characteristic-representation}
Let $\mathsf Z$ be standard normal and let
$\mathsf E_1,\mathsf E_2,\ldots$ be independent mean-one exponential
variables, independent of $\mathsf Z$.  Then
\begin{equation}
 \label{hp:eq:Xi}
 \mathsf Y_g\defeq\gamma+\sqrt q\,\mathsf Z
 +\sum_{j=1}^{\infty}y_j\left(\mathsf E_j-1\right)
\end{equation}
is well defined by $\mathrm L^2$ convergence of the series, and
\begin{equation}
 \label{hp:eq:h-characteristic}
 h_g^{\mathrm{HP}}\left(u\right)=\E\left[\mathrm e^{\mathrm i u\mathsf Y_g}\right],\qquad
 \E\left[\mathsf Y_g\right]=\gamma,\qquad
 \operatorname{Var}\left(\mathsf Y_g\right)=\delta.
\end{equation}
In particular,
\begin{equation}
 \label{hp:eq:h-basic-bounds}
 |h_g^{\mathrm{HP}}\left(u\right)|\leq1,\qquad
 h_g^{\mathrm{HP}}\left(-u\right)=\overline{h_g^{\mathrm{HP}}\left(u\right)},
\end{equation}
and, for $k=1,2$,
\begin{equation}
 \label{hp:eq:h-derivative-moment-bound}
 |\partial_u^kh_g^{\mathrm{HP}}\left(u\right)|\leq \E\left[|\mathsf Y_g|^k\right],
 \qquad
 \E\left[|\mathsf Y_g|^2\right]=\gamma^2+\delta.
\end{equation}
Moreover,
\begin{equation}
 \label{hp:eq:h-modulus}
 |h_g^{\mathrm{HP}}\left(u\right)|^2=\mathrm e^{-qu^2}\prod_{j=1}^{\infty}\left(1+y_j^2u^2\right)^{-1}
 \leq\left(1+\delta u^2\right)^{-1}.
\end{equation}
\end{lem}

\begin{proof}
The centred exponential series in \eqref{hp:eq:Xi} is an orthogonal series
in $\mathrm L^2$, with sum of variances $\sum_{j=1}^{\infty}y_j^2<\infty$.  Since
\[
 \E\left[\mathrm e^{\mathrm i u y\left(\mathsf E-1\right)}\right]
 =\frac{\mathrm e^{-\mathrm i uy}}{1-\mathrm i uy},
\]
independence and $\mathrm L^2$ convergence give
\[
 \E\left[\mathrm e^{\mathrm i u\mathsf Y_g}\right]
 =\mathrm e^{\mathrm i\gamma u-qu^2/2}
  \prod_{j=1}^{\infty}\frac{\mathrm e^{-\mathrm i uy_j}}{1-\mathrm i uy_j}
 =\frac1{g\left(\mathrm i u\right)}.
\]
This proves \eqref{hp:eq:h-characteristic}--\eqref{hp:eq:h-derivative-moment-bound}.
The identity in \eqref{hp:eq:h-modulus} follows by taking absolute values.
Finally,
\[
 \mathrm e^{qu^2}\prod_{j=1}^{\infty}\left(1+y_j^2u^2\right)
 \geq1+\left(q+\sum_{j=1}^{\infty}y_j^2\right)u^2,
\]
which proves the last inequality.
\end{proof}

Differentiating evaluation martingales recovers the state coordinates;
local covariance bounds justify this operation without global moment assumptions. This is the content of the next technical lemma.
\begin{lem}
\label{hp:lem:derivative-martingales}
Let \(\boldsymbol{\Phi}\) solve the Hua--Pickrell martingale problem.  For every
\(z_0,w_0\in\C\) and \(k,\ell\geq0\),
\begin{equation}
 \label{hp:eq:derivative-martingale}
 \partial_z^kM_t^{\boldsymbol{\Phi}}\left(z_0\right)
 =\partial_z^k\boldsymbol{\Phi}_t\left(z_0\right)-\partial_z^k\boldsymbol{\Phi}_0\left(z_0\right)
 -\int_0^t\partial_z^k
   \left(\mathcal L_{a,b}^{\mathrm{HP}}\boldsymbol{\Phi}_s\right)\left(z_0\right)\,\dd s
\end{equation}
is a continuous local martingale, and
\begin{equation}
 \label{hp:eq:derivative-bracket}
 \left[\partial_z^kM^{\boldsymbol{\Phi}}\left(z_0\right),\partial_w^\ell M^{\boldsymbol{\Phi}}\left(w_0\right)\right]_t
 =\int_0^t\partial_z^k\partial_w^\ell
 \mathcal Q_{\boldsymbol{\Phi}_s}^{\mathrm{HP}}\left(z_0,w_0\right)\,\dd s.
\end{equation}
\end{lem}

\begin{proof}
We first abbreviate for simplicity
\begin{equation*}
    M_t(z)= M_t^{\boldsymbol{\Phi}}\left(z\right),\ \mathcal{Q}_t(z,w)=\mathcal Q_{\boldsymbol{\Phi}_t}^{\mathrm{HP}}(z,w),\
    \mathcal{L}=\mathcal L_{a,b}^{\mathrm{HP}}.
\end{equation*}
Compact-open continuity of \(\boldsymbol \Phi\) and Cauchy’s estimates imply that \(t\mapsto\mathcal L\boldsymbol \Phi_t\) is compact-open continuous, hence bounded on compact space–time sets, pathwise. Cauchy’s formula and Fubini therefore give
\[
\partial_z^k\!\left(\int_0^t(\mathcal L\boldsymbol \Phi_s)(z)\,ds\right)
=
\int_0^t\partial_z^k(\mathcal L\boldsymbol \Phi_s)(z)\,ds.
\]
Thus \(M_t\) is entire in \(z\), its spatial derivatives are continuous in \(t\), and \eqref{hp:eq:derivative-martingale} holds pathwise. It remains to establish the local-martingale property and the bracket identity \eqref{hp:eq:derivative-bracket}.

Choose \(R>1+\max\{|z_0|,|w_0|\}\), and put
\[
D=\{z:|z|\le R\},\qquad D'=\{z:|z|\le R+2\}.
\]
We define the stopped process as following,
\[
\sigma_L
\defeq
L\wedge\inf\left\{
t\ge0:
\sup_{z\in D'}|\boldsymbol{\Phi}_t(z)|
+\sup_{z,w\in D'}|\mathcal{Q}_t(z,w)|
\ge L
\right\},
\qquad
M_t^L(z)=M_{t\wedge\sigma_L}(z).
\]
The quantities defining \(\sigma_L\) are adapted and continuous. Since the stopping time is common to all evaluation points, \(M_t^L\) remains entire in \(z\). Moreover, \(|Q_s(z,w)|\le L\) on \(D'\times D'\) for \(s<\sigma_L\).

Real-type symmetry and the stopped covariance bound imply
\[
[M^L(z),\overline{M^L(z)}]_T\le LT,
\qquad z\in D'.
\]
Thus each stopped evaluation is a square-integrable martingale. For any complex valued process $\mathsf X$, we define
\[
\|\mathsf X\|_{2,T}\defeq
\left(\mathbb E\left[\sup_{t\le T}|\mathsf X_t|^2\right]\right)^{1/2}.
\]
The martingale isometry and Doob’s inequality yield
\[
\mathbb E\left[|M_T^L(z)|^2\right]\le LT,
\qquad
\|M^L(z)\|_{2,T}\defeq\le2\sqrt{LT},
\]
These estimates require no initial moment assumption: if the initial supremum already reaches \(L\), then \(\sigma_L=0\) and \(M^L\equiv0\).

We next obtain spatial continuity in this norm. Cauchy’s estimates and integration of \(\partial_1\partial_2\mathcal Q_s\) along the segments \(w\to z\) and \(\overline w\to\overline z\) bound the Hermitian bracket density of \(M^L(z)-M^L(w)\) by \(L|z-w|^2\) before stopping. The martingale isometry and Doob’s inequality therefore give
\[
\|M^L(z)-M^L(w)\|_{2,T}
\le2\sqrt{LT}\,|z-w|,
\qquad z,w\in D.
\]
This estimate will allow the Cauchy contour sums to converge in \(\|\cdot\|_{2,T}\).
Put \(z_\theta=z_0+\mathrm{e}^{i\theta}\) and \(F_\theta(t)=\mathrm{e}^{-ik\theta}M^L_t(z_\theta)\). The preceding estimates imply
\[
\|F_\theta-F_\varphi\|_{2,T}
\le2(k+1)\sqrt{LT}\,|\theta-\varphi|.
\]
By Cauchy’s formula,
\[
A_t\defeq \partial_z^kM_t^L(z_0)
=\frac{k!}{2\pi}\int_0^{2\pi}F_\theta(t)\,d\theta.
\]
For \(\delta_N=2\pi/N\) and \(\theta_j=j\delta_N\), we have
\[
A_t^{(N)}
\defeq
\frac{k!}{N}\sum_{j=0}^{N-1}
F_{\theta_j}(t),\quad \sup_{N\in\N}\|A^{(N)}\|_{2,T}\le2k!\sqrt{LT}.
\]
Each \(A^{(N)}\) is a continuous square-integrable martingale. Minkowski’s inequality and the preceding Lipschitz estimate give
\[
\begin{aligned}
\|A^{(N)}-A\|_{2,T}
\le
\frac{k!}{2\pi}
\sum_{j=0}^{N-1}
\int_{\theta_j}^{\theta_{j+1}}
\|F_{\theta_j}-F_\theta\|_{2,T}\,d\theta
\le2k!(k+1)\sqrt{LT}\,\delta_N
\longrightarrow0.
\end{aligned}
\]
Thus \(A\) is adapted and square-integrable. Let $||\cdot||_2$ be the $\mathrm L^2$ norm. For \(0\le s\le t\le T\), the martingale property of \(A^{(N)}\) and the \(\mathrm L^2\)-contraction of conditional expectation yield
\[
\left\|\mathbb E[A_t\mid\mathscr F_s]-A_s\right\|_2
\le
\|A_t-A_t^{(N)}\|_2
+\|A_s^{(N)}-A_s\|_2
\le2\|A-A^{(N)}\|_{2,T}
\longrightarrow0.
\]
Hence \(A\) is a continuous square-integrable martingale. Similarly, $B_t\defeq\partial_w^\ell M_t^L(w_0)$ is a continuous square-integrable martingale.

Set \(C^{(N)}=[A^{(N)},B^{(N)}]\). By the bilinear covariance identity, \(C^{(N)}\) is the corresponding double-contour Riemann sum involving \(\mathcal Q_s\). Cauchy’s formula and joint continuity identify its limit as
\[
C_t
=
\int_0^{t\wedge\sigma_L}
\partial_z^k\partial_w^\ell
\mathcal Q_s(z_0,w_0)\,ds,
\]
with convergence uniform on bounded time intervals, almost surely. The localized covariance bound gives
\(\sup_{t\le T}|C_t^{(N)}|\le k!\ell!LT\), so dominated convergence yields
\[
\mathbb E\left[\sup_{t\le T}|C_t^{(N)}-C_t|\right]\longrightarrow0.
\]
Similarly, the preceding \(\|\cdot\|_{2,T}\)-convergence and Cauchy–Schwarz imply
\[
\E \left[\sup_{t\le T}
|A_t^{(N)}B_t^{(N)}-A_tB_t|\right]
\longrightarrow0.
\]

By Itô’s product formula, \(A^{(N)}B^{(N)}-C^{(N)}\) is a local martingale. The preceding bounds give it an integrable supremum on every finite interval, hence it is a true martingale. Passing its conditional-expectation identities to the limit in \(\mathrm L^1\) shows that \(AB-C\) is a martingale. Since \(AB-[A,B]\) is also a local martingale, \([A,B]-C\) is a continuous local martingale of finite variation, starting at zero. It therefore vanishes, proving
\[
[A,B]_t
=
\int_0^{t\wedge\sigma_L}
\partial_z^k\partial_w^\ell
\mathcal Q_s(z_0,w_0)\,ds.
\]

Finally, compact-open path continuity and continuity of the covariance kernel imply that the quantities defining \(\sigma_L\) are bounded on every finite time interval, pathwise. Hence \(\sigma_L\uparrow\infty\). Because the stopping time is independent of the spatial variables, \(A\) and \(B\) are the stopped versions of the corresponding derivatives of \(M^{\boldsymbol\Phi}\). These derivatives are therefore continuous local martingales, and removing the stopping in the preceding identity gives \eqref{hp:eq:derivative-bracket}. Together with the pathwise identity \eqref{hp:eq:derivative-martingale}, this proves the lemma. \(\square\)
\end{proof}
To estimate the martingale part of the coordinate process in the next lemma, we first record a factorisation of the Hua–Pickrell covariance. Retain the measurable signed enumeration
\eqref{hp:eq:signed-root-enumeration}.  Logarithmic differentiation of
\eqref{eq:intro-HP-state} gives
\begin{equation}
 \label{hp:eq:log-derivative}
 \frac{g'\left(z\right)}{g\left(z\right)}=-\gamma-qz-
 \sum_{j=1}^{\infty}\frac{y_j^2z}{1-y_jz}.
\end{equation}
Consequently,
\begin{equation}
 \label{hp:eq:cov-factorisation}
 2\frac{zw}{z-w}\left(g\left(z\right)g'\left(w\right)-g\left(w\right)g'\left(z\right)\right)
 =2zwg\left(z\right)g\left(w\right)\left(q+
 \sum_{j=1}^{\infty}\frac{y_j^2}{\left(1-y_jz\right)\left(1-y_jw\right)}\right).
\end{equation}
The following algebraic square root of the covariance will be used below:
\begin{equation}
 \label{hp:eq:noise-coefficients}
 b_0\left(g;z\right)\defeq-\sqrt{2q}\,zg\left(z\right),
 \qquad
 b_j\left(g;z\right)\defeq-\sqrt2\,\frac{zy_jg\left(z\right)}{1-y_jz},\quad j\geq1.
\end{equation}
The apparent poles are removable, and \eqref{hp:eq:cov-factorisation}
shows that the products of the coefficients in
\eqref{hp:eq:noise-coefficients} sum to the prescribed covariance.

Set $R(g)=\gamma(g)^2+\delta(g)$.

The first two coordinates satisfy a closed growth estimate that controls
probe derivatives, including after truncation on an initial event.
\begin{lem}
\label{hp:lem:state-moment}
Let \(\boldsymbol{\Phi}\) be a continuous weak solution of the
Hua--Pickrell martingale problem
\eqref{eq:intro-martingale}--\eqref{eq:intro-bracket}
with parameters \(a,b\in\R\). Write
\[
 \mathsf\Gamma_t\defeq\gamma\left(\boldsymbol{\Phi}_t\right),
 \qquad
 \mathsf\Delta_t\defeq\delta\left(\boldsymbol{\Phi}_t\right),
 \qquad
 \mathsf R_t\defeq R\left(\boldsymbol{\Phi}_t\right).
\]
Then
\begin{align}
 \dd\mathsf\Gamma_t&=\dd\mathsf N_t +2\left(b-a\mathsf\Gamma_t\right)\,\dd t, \qquad \dd\left[\mathsf N\right]_t=2\mathsf\Delta_t\,\dd t,
 \label{hp:eq:gamma-dynamics}\\
 \operatorname{drift}\left(\dd\mathsf\Delta_t\right)&= \left(2+2\mathsf\Gamma_t^2-4a\mathsf\Delta_t\right)\,\dd t,\qquad \operatorname{drift}\left(\dd\mathsf R_t\right)= \left(2+4b\mathsf\Gamma_t+\left(2-4a\right)\mathsf R_t\right)\,\dd t.
 \label{hp:eq:delta-drift}
\end{align}
The martingale part $\mathsf N^R$ of $\mathsf R$ satisfies
\begin{equation}
 \label{hp:eq:R-bracket}
 \dd\left[\mathsf N^R\right]_t\leq16\mathsf R_t^2\,\dd t.
\end{equation}
Consequently, for every $m\geq1$, $T<\infty$, and
${\mathscr A}\in\mathscr F_0$ such that
$\E\left[\1_{\mathscr A}(1+\mathsf R_0)^m\right]<\infty$,
\begin{equation}
 \label{hp:eq:R-moment}
 \E\left[\1_{\mathscr A}\sup_{t\leq T}\left(1+\mathsf R_t\right)^m\right]
 \leq C_{m,T,a,b}\E\left[\1_{\mathscr A}\left(1+\mathsf R_0\right)^m\right].
\end{equation}
\end{lem}

\begin{proof}
Lemma~\ref{hp:lem:derivative-martingales} allows differentiation of the
martingale identities.  The first two derivatives of the drift at zero are
\[
 -2b-2ag'\left(0\right),\qquad
 4b\mathsf\Gamma-2-2\left(2a+1\right)\boldsymbol{\Phi}''\left(0\right),
\]
and the first mixed derivative of the bracket is \(2\delta\).  These give
\eqref{hp:eq:gamma-dynamics}; It\^o's formula for
\(\mathsf\Delta=\mathsf\Gamma^2-\boldsymbol{\Phi}''(0)\) and
\(\mathsf R=\mathsf\Gamma^2+\mathsf\Delta\) gives
\eqref{hp:eq:delta-drift}.
The covariance factorisation \eqref{hp:eq:noise-coefficients} gives
\[
 \frac{\dd\left[\mathsf N^R\right]_t}{\dd t}
 =8q\mathsf\Gamma^2+8\sum_{j=1}^{\infty}\left(\mathsf\Gamma y_j+y_j^2\right)^2
 \leq16\left(\mathsf\Gamma^2+\mathsf\Delta\right)^2,
\]
using \(\sum_{j=1}^{\infty}y_j^4\leq\mathsf\Delta^2\).

Put \(\mathsf X=1+\mathsf R\) and
\(\tau_L=\inf\{t:\mathsf R_t\geq L\}\).  Its drift is bounded in absolute
value by \(C_{a,b}\mathsf X\), and
\(\dd[\mathsf N^R]\leq16\mathsf X^2\dd t\).  It\^o's formula for
\(\mathsf X^m\), followed by the Burkholder--Davis--Gundy and Young
inequalities, gives
\[
 \E\left[\sup_{s\leq t}\left|\1_{\mathscr A} m\int_0^{s\wedge\tau_L} \mathsf X_v^{m-1}\,\dd\mathsf N_v^R\right|\right] \leq\frac12\E\left[\1_{\mathscr A}\sup_{s\leq t} \mathsf X_{s\wedge\tau_L}^m\right] +C_m\int_0^t\E\left[\1_{\mathscr A}\mathsf X_{v\wedge\tau_L}^m\right]\,\dd v.
\]
Here \({\mathscr A}\in\mathscr F_0\) permits multiplication of the martingale
integrand by \(\1_{\mathscr A}\).  Absorbing the supremum term and applying Gronwall
gives \eqref{hp:eq:R-moment} up to \(\tau_L\), uniformly in \(L\).
Path continuity gives \(\tau_L\uparrow\infty\), and Fatou removes the stopping.
\end{proof}
\subsection{Positive bounded probes}
Similarly to the Hua--Pickrell case, the reciprocals on the negative axis form bounded observables in the positive case.
\begin{defn}\label{def probe h positive}
    For $g\in\mathscr E_+$, and $s>0$, define
\begin{equation*}
h_g^{+}\left(s\right)\defeq 1/g\left(\mathrm -s\right),\qquad s>0.
\end{equation*}
For fixed \(s\), this is continuous in \(g\) for locally uniform
convergence.
\end{defn}
We also have the following lemma.
\begin{lem}
    If $(\mathsf E_j)$ are independent mean-one exponential variables, then
\begin{equation}
 \label{eq:Laplace-representation}
 h_g^+\left(s\right)=\E\left[\mathrm e^{-s\mathsf Y_g^+}\right],\qquad
 \mathsf Y_g^+\defeq d\left(g\right)+\sum_{j=1}^{\infty}x_j\mathsf E_j.
\end{equation}
Moreover, for all $m\geq0$, we have
\begin{equation}
 \label{eq:h-derivative-bound}
 0<h_g^+\left(s\right)\leq1,\qquad
 \sup_{g\in\mathscr E_+}|\partial_s^mh_g^+\left(s\right)|
 \leq\sup_{y\geq0}y^m\mathrm e^{-sy}
 =\left(\frac{m}{\mathrm e s}\right)^m,
\end{equation}
with value $1$ when $m=0$.  
\end{lem}
\begin{proof}
    Indeed, the nonnegative random variable $\mathsf Y_g^+$ has finite
expectation $\gamma(g)$, and independence gives the formula first for
partial products and then by bounded convergence in general. For all $m\geq0$, differentiation under
the expectation gives \eqref{eq:h-derivative-bound}, hence completing the proof.
\end{proof}

\section{Reciprocal dualities and weak uniqueness}
\label{sec:duality-uniqueness}

\subsection{Positive dual}

The product of reciprocals on the negative axis determine the positive state law. Letting $h_g^+$ be the probe defined in Definition \ref{def probe h positive}, put
\begin{equation}
 \label{eq:reciprocal-probes}
 H_n^+\left(g,\mathbf s\right)\defeq\prod_{i=1}^{n}h_g^+\left(s_i\right), \quad \mathbf{s}\in(0,\infty)^n.
\end{equation}

On $\mathscr W_n=\{0<s_1<\cdots<s_n\}$, define
\begin{equation}
 \mathcal G_n^\theta \defeq\frac12\sum_{i=1}^{n} s_i^2\partial_{ii} +\frac\theta2\sum_{i=1}^{n} s_i\partial_i +\sum_{1\leq i<j\leq n}\frac{s_is_j}{s_i-s_j}\left(\partial_i-\partial_j\right),\qquad V_n^+\left(\mathbf s\right)\defeq\frac12\sum_{i=1}^{n} s_i.
 \label{eq:dual-generator}
\end{equation}

The positive dual evolves finitely many probe arguments \(\mathbf s\).
\begin{lem}
\label{lem:finite-dual}
For every $\theta\in\R$, $\mathcal G_n^\theta$ generates a conservative
continuous diffusion $\boldsymbol{\mathsf S}=(\mathsf S_i)_{i=1}^n$ in $\mathscr W_n$, satisfying
\begin{equation}
 \label{eq:dual-SDE}
 \dd\mathsf S_i=\mathsf S_i\,\dd\mathsf B_i+
 \left(\frac\theta2\mathsf S_i+
 \sum_{\substack{j=1\\j\neq i}}^{n}\frac{\mathsf S_i\mathsf S_j}
                         {\mathsf S_i-\mathsf S_j}\right)\dd t.
\end{equation}
Moreover, for \(m\ge1\), define
\[
\rho_m\defeq\inf\{t\ge0:\boldsymbol{\mathsf S}(t)\notin \mathscr W_n
                    \text{ or }\mathsf S_n(t)\ge m\},
\qquad
\rho\defeq\lim_{m\to\infty}\rho_m,
\]
with \(\inf\varnothing=\infty\).
Then \(\mathbb P_s(\rho=\infty)=1\) for every \(s\in \mathscr W_n\).
\end{lem}

\begin{proof} We will do the construction in logarithmic coordinates. Put $y_i=\log s_i$, $\mu=(\theta-n)/2$, and
\[
 \mathscr{V}_n\left(\mathbf y\right)\defeq\prod_{1\leq i<j\leq n}2\sinh\left(\frac{y_j-y_i}{2}\right)
 =\det\left[\mathrm e^{\eta_i y_j}\right]_{i,j=1}^n,
 \qquad \eta_i=i-\frac{n+1}{2}.
\]
The identity $s_j/(s_i-s_j)=\frac12\coth((y_i-y_j)/2)-\frac12$
and It\^o's formula give the target generator in logarithmic coordinates
\begin{equation}
 \label{eq:log-dual-generator}
 \tilde{\mathcal G}_n^\theta=\frac12\Delta_{\mathbf y}
 +\mu\sum_{i=1}^{n}\partial_{y_i}+\nabla\log\mathscr{V}_n\cdot\nabla.
\end{equation}
Every determinant monomial has the exponents $\eta_i$ in some order, so
\[
 \frac12\Delta_{\mathbf y}\mathscr{V}_n=\varkappa_n\mathscr{V}_n,
 \qquad\sum_{i=1}^{n}\partial_{y_i}\mathscr{V}_n=0,
 \qquad\varkappa_n=\frac12\sum_{i=1}^{n}\eta_i^2=\frac{n\left(n^2-1\right)}{24}.
\]
Let $\boldsymbol{\mathsf Z}$ be independent Brownian motions with common
drift $\mu$, and define the first collision time
\begin{equation*}
    \tau_\partial\defeq\inf\{t\ge0:\mathsf Z_i(t)=\mathsf Z_j(t)\text{ for some }1\le i<j\le n\}
\end{equation*}
Karlin--McGregor \cite{pjm/1103038889} and Andr\'eief \cite{Forrester2019} identities give
\[
 \E_{\mathbf y}\left[\mathscr{V}_n\left(\boldsymbol{\mathsf Z}_t\right);t<\tau_\partial\right]
 =\det\left[\mathrm e^{\eta_i y_j+\mu\eta_i t+\eta_i^2t/2}\right]_{i,j=1}^n
 =\mathrm e^{\varkappa_nt}\mathscr{V}_n\left(\mathbf y\right).
\]
Here the product of the two determinants is symmetric, so the factor
$1/n!$ from chamber integration cancels the factor in Andr\'eief's
identity; all integrals are finite Gaussian exponential moments.
Consequently,
\[
 \mathbf Q_tf\left(\mathbf y\right)=\frac{\mathrm e^{-\varkappa_nt}}{\mathscr{V}_n\left(\mathbf y\right)}
 \E_{\mathbf y}\left[\mathscr{V}_n\left(\boldsymbol{\mathsf Z}_t\right)
 f\left(\boldsymbol{\mathsf Z}_t\right);t<\tau_\partial\right]
\]
is a probability semigroup. Its conservative Doob $h$-transform
\cite{DoobPotential,RevuzYor} has
continuous paths staying in the open chamber at every finite time and
generator \eqref{eq:log-dual-generator}.  Exponentiation proves the claim.
\end{proof}

\subsection{The Hua--Pickrell signed dual}

Using the probes from Lemma~\ref{hp:lem:characteristic-representation}, put
\begin{align}
 \label{hp:eq:dual-notation}
 H_n^{\mathrm{HP}}\left(g,\mathbf u\right)&\defeq\prod_{i=1}^{n}h_g^{\mathrm{HP}}\left(u_i\right),
 &V_n^{\mathrm{HP}}\left(\mathbf u\right)&\defeq\sum_{i=1}^{n} u_i^2,\\
 C_{b,n}\left(\mathbf u\right)&\defeq2\mathrm i b\sum_{i=1}^{n} u_i-V_n^{\mathrm{HP}}\left(\mathbf u\right).
 \nonumber
\end{align}
For pairwise distinct nonzero coordinates define
\begin{equation}
 \label{hp:eq:dual-generator}
 \mathcal A_n^a\defeq\sum_{i=1}^{n} u_i^2\partial_{ii}
 -2a\sum_{i=1}^{n} u_i\partial_i
 +2\sum_{1\leq i<j\leq n}\frac{u_iu_j}{u_i-u_j}\left(\partial_i-\partial_j\right).
\end{equation}
The signed chambers are
\begin{equation}
 \label{hp:eq:signed-chamber}
 \mathscr C_{\ell,n}\defeq
 \{u_1<\cdots<u_\ell<0<u_{\ell+1}<\cdots<u_n\},
 \qquad0\leq\ell\leq n.
\end{equation}
Any vector \(\mathbf{u}\in\mathbb R^n\) with pairwise distinct nonzero coordinates has a unique increasing rearrangement in \(\mathscr C_{\ell,n}\), where \(\ell=\#\{i:u_i<0\}\). All the duality formulas below are invariant under this relabelling.

The signed dual preserves each argument's sign and avoids collisions,
providing the corresponding Feynman--Kac representation for Hua--Pickrell probes.
\begin{lem}
\label{hp:lem:signed-dual}
On every $\mathscr C_{\ell,n}$, $\mathcal A_n^a$ generates a conservative
continuous diffusion satisfying
\begin{equation}
 \label{hp:eq:dual-SDE}
 \dd\mathsf S_i=\sqrt2\mathsf S_i\,\dd\mathsf B_i+
 \left[-2a\mathsf S_i+2\sum_{\substack{j=1\\j\neq i}}^{n}
 \frac{\mathsf S_i\mathsf S_j}{\mathsf S_i-\mathsf S_j}\right]\dd t.
\end{equation}
Moreover, for \(m\ge1\), define
\[
\rho_m\defeq\inf\left\{t\ge0:\boldsymbol{\mathsf S}(t)\notin\mathscr C_{\ell,n}
                    \text{ or }\max_{1\le i\le n}|\mathsf S_i(t)|\ge m\right\},
\qquad
\rho\defeq\lim_{m\to\infty}\rho_m,
\]
with \(\inf\varnothing=\infty\).
Then \(\mathbb P_u(\rho=\infty)=1\) for every
\(u\in\mathscr C_{\ell,n}\).
\end{lem}

\begin{proof}
Set
\[
 \mathscr{V}\left(\mathbf u\right)\defeq\det\left[u_i^{j-1}\right]_{i,j=1}^n
 =\prod_{1\leq i<j\leq n}\left(u_j-u_i\right),\quad
 \beta\defeq-2a-2\left(n-1\right),\quad
 \mathcal L_0\defeq\sum_{i=1}^{n}\left(u_i^2\partial_{ii}+\beta u_i\partial_i\right),
\]
and $\Lambda_n=\sum_{r=0}^{n-1}r(r-1+\beta)$.  The monomial
expansion gives $\mathcal L_0\mathscr V =\Lambda_n\mathscr V$; the product rule,
$\partial_i\log\mathscr V=\sum_{\substack{j=1\\j\neq i}}^{n}(u_i-u_j)^{-1}$, and
$u_i^2/(u_i-u_j)=u_i+u_iu_j/(u_i-u_j)$ give
\begin{equation}
 \label{hp:eq:dual-Doob}
 \mathcal A_n^af=\mathscr V^{-1}\mathcal L_0\left(\mathscr V f\right)-\Lambda_nf.
\end{equation}
Let $p_t^\beta(x,y)$ be the transition density of signed geometric
Brownian motion with generator $x^2\partial_{xx}+\beta x\partial_x$ on
$\R^*=\R\setminus\{0\}$.  It vanishes when $xy<0$; hence
Karlin--McGregor \cite{pjm/1103038889} applied separately to the negative and positive blocks
gives the killed density $\det\left[p_t^\beta(u_i,v_j)\right]$ on
$\mathscr C_{\ell,n}$.  Its product with $\mathscr V(\mathbf v)$ is symmetric
in $\mathbf v$ and vanishes unless precisely $\ell$ terminal coordinates
are negative.  That labelled region is the union of the $n!$ permutations
of the chamber.  The moment identity
\[
 \int_{\R^*}p_t^\beta\left(u,v\right)v^r\,\dd v
 =\mathrm e^{r\left(r-1+\beta\right)t}u^r
\]
and Andr\'eief identity \cite{Forrester2019} therefore give
\begin{equation}
 \int_{\mathscr C_{\ell,n}} \det\left[p_t^\beta\left(u_i,v_j\right)\right]\mathscr V\left(\mathbf v\right)\,\dd\mathbf v =\det\left[\int_{\R^*} p_t^\beta\left(u_i,v\right)v^{j-1}\,\dd v\right]_{i,j=1}^n =\mathrm e^{\Lambda_nt}\mathscr V\left(\mathbf u\right).
 \label{hp:eq:dual-conservative-integral}
\end{equation}
The integrals are absolutely convergent because geometric Brownian motion
has every finite absolute moment.  Since $\mathscr V>0$ in the chamber,
\eqref{hp:eq:dual-Doob}--\eqref{hp:eq:dual-conservative-integral} define a
conservative Doob \(h\)-transform.  Its continuous paths remain in the chamber
at every finite time, and its generator gives the displayed SDE.
\end{proof}

\subsection{Conditional duality}

The following identities express conditional probe expectations through
the finite duals, without assuming a Markov property for the original solution. This is the heart of the argument.
\begin{prop}
\label{prop:conditional-duality}
For every initial law and $0\leq r\leq t$, every solution of the positive
$(\theta,\kappa)$ martingale problem satisfies
\begin{equation}
 \label{eq:conditional-duality}
 \E\left[H_n^+\left(\boldsymbol{\Phi}_t,\mathbf s\right)\mid\mathscr F_r\right]
 =\E_{\mathbf s}\left[
 \mathrm e^{-\kappa\int_0^{t-r}V_n^+\left(\boldsymbol{\mathsf S}_v\right)\,\dd v}H_n^+\left(\boldsymbol{\Phi}_r,\boldsymbol{\mathsf S}_{t-r}\right)
 \right],
\end{equation}
where the expectation on the right is taken only with respect to the (independent) dual diffusion with generator $\mathcal{G}_n^\theta$.
Furthermore, every solution of the Hua--Pickrell martingale problem satisfies
\begin{equation}
 \label{hp:eq:conditional-duality}
 \E\left[H_n^{\mathrm{HP}}\left(\boldsymbol{\Phi}_t,\mathbf u\right)\mid\mathscr F_r\right]
 =\E_{\mathbf u}\left[
 \mathrm e^{\int_0^{t-r}C_{b,n}\left(\boldsymbol{\mathsf S}_v\right)\,\dd v}
 H_n^{\mathrm{HP}}\left(\boldsymbol{\Phi}_r,\boldsymbol{\mathsf S}_{t-r}\right)\right].
\end{equation}
for pairwise distinct nonzero $\mathbf u$, where the expectation on the
right is over the (independent) signed dual diffusion with generator $\mathcal{A}_n^a$.
\end{prop}

\begin{proof}
Write $\mathsf r_t(z)=\boldsymbol{\Phi}_t(z)^{-1}$ on
$\C\setminus[0,\infty)$ in the positive case and on $\C\setminus\R$ in
the Hua--Pickrell case.  These are fixed zero-free domains.  For each
evaluation, continuity and nonvanishing make the stopping times before
the denominator approaches zero tend to infinity.  Holomorphic It\^o
calculus, with $\dd N_t^{\boldsymbol{\Phi}}(z)=-\boldsymbol{\Phi}_t(z)^{-2}\dd M_t^{\boldsymbol{\Phi}}(z)$, gives,
in the positive and Hua--Pickrell cases respectively,
\begin{align}
 \label{eq:reciprocal-Ito}
 \dd\mathsf r_t\left(z\right)&=\dd N_t^{\boldsymbol{\Phi}}\left(z\right)+
 \left(\frac\kappa2z\mathsf r_t\left(z\right)+\frac\theta2z\mathsf r_t'\left(z\right)
       +\frac12z^2\mathsf r_t''\left(z\right)\right)\dd t,\\
 \label{hp:eq:reciprocal-Ito}
 \dd\mathsf r_t\left(z\right)&=\dd N_t^{\boldsymbol{\Phi}}\left(z\right)+
 \left(\left(2bz+z^2\right)\mathsf r_t\left(z\right)-2az\mathsf r_t'\left(z\right)
       +z^2\mathsf r_t''\left(z\right)\right)\dd t.
\end{align}
Indeed, $\mathsf r'=-\boldsymbol{\Phi}'/\boldsymbol{\Phi}^2$ and
$\mathsf r''=2(\boldsymbol{\Phi}')^2/\boldsymbol{\Phi}^3-\boldsymbol{\Phi}''/\boldsymbol{\Phi}^2$;
the diagonal bracket supplies the second derivative term in each formula.
For distinct arguments, the cross-bracket is
\begin{equation}
 \label{eq:reciprocal-cross-bracket}
 \dd\left[N^{\boldsymbol{\Phi}}\left(z\right),N^{\boldsymbol{\Phi}}\left(w\right)\right]_t
 =c_\diamond\frac{zw}{z-w}
 \left(\mathsf r_t'\left(z\right)\mathsf r_t\left(w\right)
       -\mathsf r_t\left(z\right)\mathsf r_t'\left(w\right)\right)\,\dd t,
 \qquad c_+=1,\quad c_{\mathrm{HP}}=2.
\end{equation}
Put
\[
 h_t^+(s)\defeq h_{\boldsymbol{\Phi}_t}^+(s)
 =\mathsf r_t(-s),\qquad
 h_t^{\mathrm{HP}}(u)
 \defeq h_{\boldsymbol{\Phi}_t}^{\mathrm{HP}}(u)
 =\mathsf r_t(\mathrm i u),
\]
and let
\[
 L_t^+(s)\defeq N_t^{\boldsymbol{\Phi}}(-s),\qquad
 L_t^{\mathrm{HP}}(u)
 \defeq N_t^{\boldsymbol{\Phi}}(\mathrm i u).
\]
The chain rule gives
\[
 \begin{aligned}
 (h_t^+)'(s)&=-\mathsf r_t'(-s),
 & (h_t^+)''(s)&=\mathsf r_t''(-s),\\
 (h_t^{\mathrm{HP}})'(u)&=\mathrm i\mathsf r_t'(\mathrm i u),
 & (h_t^{\mathrm{HP}})''(u)&=-\mathsf r_t''(\mathrm i u).
 \end{aligned}
\]
Thus the two reciprocal It\^o identities become
\begin{align*}
 \dd h_t^+(s)
 &=\dd L_t^+(s)
   +\left(\frac12s^2\partial_s^2
          +\frac\theta2s\partial_s-\frac\kappa2s\right)
     h_t^+(s)\,\dd t,\\
 \dd h_t^{\mathrm{HP}}(u)
 &=\dd L_t^{\mathrm{HP}}(u)
   +\left(u^2\partial_u^2-2au\partial_u
          +2\mathrm i bu-u^2\right)
     h_t^{\mathrm{HP}}(u)\,\dd t.
\end{align*}
Likewise, substituting into
\eqref{eq:reciprocal-cross-bracket} yields, for distinct arguments,
\[
 \dd[L^\diamond(v_i),L^\diamond(v_j)]_t
 =c_\diamond\frac{v_iv_j}{v_i-v_j}
 \left((h_t^\diamond)'(v_i)h_t^\diamond(v_j)
       -h_t^\diamond(v_i)(h_t^\diamond)'(v_j)\right)\dd t.
\]
Here primes on \(h_t^\diamond\) denote derivatives in its real
argument, and \(c_+=1\), \(c_{\mathrm{HP}}=2\).

Fix \(\mathbf v\) in the corresponding chamber and put
\[
 H_t^\diamond(\mathbf v)
 \defeq H_n^\diamond(\boldsymbol{\Phi}_t,\mathbf v)
 =\prod_{i=1}^{n}h_t^\diamond(v_i).
\]
The product rule for continuous semimartingales gives
\begin{align*}
 \dd H_t^\diamond(\mathbf v)
 ={}&\sum_{i=1}^{n}
       \left(\prod_{\substack{k=1\\k\neq i}}^{n}
                    h_t^\diamond(v_k)\right)
       \dd h_t^\diamond(v_i)+\sum_{1\leq i<j\leq n}
       \left(\prod_{\substack{k=1\\k\notin\{i,j\}}}^{n}
                    h_t^\diamond(v_k)\right)
       \dd[L^\diamond(v_i),L^\diamond(v_j)]_t.
\end{align*}
Empty products are understood to equal one. Since
\[
 \begin{aligned}
 \partial_{ii}H_t^\diamond
 &=(h_t^\diamond)''(v_i)
   \prod_{\substack{k=1\\k\neq i}}^{n}h_t^\diamond(v_k),\\
 (\partial_i-\partial_j)H_t^\diamond
 &=\left((h_t^\diamond)'(v_i)h_t^\diamond(v_j)
        -h_t^\diamond(v_i)(h_t^\diamond)'(v_j)\right)
   \prod_{\substack{k=1\\k\notin\{i,j\}}}^{n}
           h_t^\diamond(v_k),
 \end{aligned}
\]
the pairwise covariations give exactly the interaction terms
in the dual generators. The remaining drift terms give the
one-coordinate operators and the multiplication terms. Hence, with
\[
 (\mathcal B_+,C_+)=(\mathcal G_n^\theta,-\kappa V_n^+),
 \qquad
 (\mathcal B_{\mathrm{HP}},C_{\mathrm{HP}})
   =(\mathcal A_n^a,C_{b,n}),
\]
the process
\begin{equation}
 \label{eq:probe-martingale}
 \mathsf Z_t^\diamond\defeq
 H_n^\diamond(\boldsymbol{\Phi}_t,\mathbf v)
 -H_n^\diamond(\boldsymbol{\Phi}_0,\mathbf v)
 -\int_0^t(\mathcal B_\diamond+C_\diamond)
      H_n^\diamond(\boldsymbol{\Phi}_s,\cdot)(\mathbf v)\,\dd s
\end{equation}
has the representation
\[
 \mathsf Z_t^\diamond
 =\sum_{i=1}^{n}\int_0^t
    \left(\prod_{\substack{k=1\\k\neq i}}^{n}
                 h_s^\diamond(v_k)\right)
    \dd L_s^\diamond(v_i).
\]
The integrands are continuous, adapted and bounded in modulus by one,
so this is a continuous local martingale for every fixed
\(\mathbf v\) in the relevant chamber.

For the positive probes, \eqref{eq:h-derivative-bound} bounds all spatial
derivatives uniformly over $g$ on every compact subchamber.  Thus
$\mathsf Z^+$ is a true martingale on bounded time intervals.
For Hua--Pickrell probes, Lemma~\ref{hp:lem:characteristic-representation}
gives
\[
 |\partial_iH_n^{\mathrm{HP}}|\leq R^{1/2},\qquad
 |\partial_{ii}H_n^{\mathrm{HP}}|\leq R,\qquad
 |\partial_{ij}H_n^{\mathrm{HP}}|\leq R\quad\left(i\ne j\right),
\]
and hence, for every compact subchamber ${\mathscr D}$,
\begin{equation}
 \label{hp:eq:probe-generator-bound}
 \sup_{\mathbf u\in {\mathscr D}}
 |\left(\mathcal A_n^a+C_{b,n}\right)H_n^{\mathrm{HP}}\left(g,\cdot\right)\left(\mathbf u\right)|
 \leq C_{\mathscr D}\left(1+R\left(g\right)\right).
\end{equation}
If ${\mathscr B}\in\mathscr F_0$ and $\E\left[\1_{\mathscr B}(1+R(\boldsymbol{\Phi}_0))^2\right]<\infty$,
Lemma~\ref{hp:lem:state-moment} shows that $\1_{\mathscr B}\mathsf Z^{\mathrm{HP}}$
is a true martingale: stop at $\tau_L=\inf\{t:R(\boldsymbol{\Phi}_t)\geq L\}$,
where the stopped process is bounded on each finite time interval, and
then use the integrable bound
\[
 \sup_{t\leq T}|\1_{\mathscr B}\mathsf Z_{t\wedge\tau_L}^{\mathrm{HP}}|
 \leq\1_{\mathscr B}\left(2+C_{\mathscr D}T\sup_{s\leq T}\left(1+R\left(\boldsymbol{\Phi}_s\right)\right)\right)
\]
to remove the stopping by dominated convergence.  In particular,
$\mathsf Z^{\mathrm{HP}}$ is a true martingale when the unweighted initial
second moment is finite.

We now prove both conditional identities at once.  Take ${\mathscr B}$ to be the whole
sample space in the positive case and ${\mathscr B}={\mathscr B}_K=\{R(\boldsymbol{\Phi}_0)\leq K\}$ in the
Hua--Pickrell case.  For ${\mathscr A}\in\mathscr F_r$ and $0\leq q\leq T$, put
\[
 U\left(q,\mathbf v\right)\defeq
 \E\left[\1_{\mathscr A}\1_{\mathscr B}H_n^\diamond\left(\boldsymbol{\Phi}_{r+q},\mathbf v\right)\right].
\]
The first two spatial derivatives of the probes are continuous in
$(g,\mathbf v)$ by locally uniform convergence and Cauchy's formula.
The bounds above, and the event-weighted moment bound of
Lemma~\ref{hp:lem:state-moment} in the Hua--Pickrell case, justify
differentiation under the expectation and continuity of these derivatives
in $(q,\mathbf v)$ on compact subchambers.  Taking expectations of the increment of the
true-martingale $\mathsf Z_t^\diamond$ between $r$ and $r+q$ gives
\begin{equation}
 \label{eq:conditional-duality-integral-PDE}
 U\left(q,\mathbf v\right)=U\left(0,\mathbf v\right)+
 \int_0^q\left(\mathcal B_\diamond+C_\diamond\right)U\left(v,\mathbf v\right)\,\dd v.
\end{equation}
The integrand is continuous, so
\begin{equation}
 \label{hp:eq:probe-backward-PDE}
 U\in C^{1,2},\qquad
 \partial_qU=\left(\mathcal B_\diamond+C_\diamond\right)U,
 \qquad |U|\leq\mathbb P\left({\mathscr A}\cap {\mathscr B}\right),
\end{equation}
locally up to the time endpoints.  For fixed $q$, backward It\^o calculus
applied to
\[
 \mathrm e^{\int_0^v C_\diamond\left(\boldsymbol{\mathsf S}_s\right)\,\dd s}
 U\left(q-v,\boldsymbol{\mathsf S}_v\right),\qquad0\leq v\leq q,
\]
has zero drift.  The displayed process is therefore a continuous local martingale on $[0,q]$. By \eqref{hp:eq:probe-backward-PDE} and $\Re C_\diamond\le0$, its modulus is bounded by $\mathbb P\left({\mathscr A}\cap {\mathscr B}\right)=1$, so it is a true martingale. Taking expectations at time $q$ gives
\[
 U\left(q,\mathbf v\right)=\E_{\mathbf v}\left[
 \mathrm e^{\int_0^q C_\diamond\left(\boldsymbol{\mathsf S}_s\right)\,\dd s}
 U\left(0,\boldsymbol{\mathsf S}_q\right)\right].
\]
Fubini gives the desired identity tested against $\1_{\mathscr A}\1_{\mathscr B}$.
In the positive case this removes ${\mathscr B}$ immediately.  For Hua--Pickrell,
$\1_{{\mathscr B}_K}\uparrow1$ and both tested integrands are bounded by one, so bounded convergence
removes ${\mathscr B}_K$.  Arbitrariness of ${\mathscr A}\in\mathscr F_r$ proves both claims
for every initial law.
\end{proof}

\begin{rmk} [On the choice and use of probes]
The form of the reciprocal probes can be understood by considering,
in the positive case, observables of the more general form
$\prod_{i=1}^n f(\boldsymbol{\Phi}_t(-s_i))$.
For $n=1$, write
\[
 G_t(s)=\boldsymbol{\Phi}_t(-s),
 \qquad F_t(s)=f(G_t(s)),
 \qquad s>0,
\]
where $f\in C^2((0,\infty))$.
Using
$\mathcal Q_g(z,z)=z^2(g'(z)^2-g(z)g''(z))$,
the drift and covariance prescribed by the martingale problem give,
by It\^o's formula,
\[
\begin{aligned}
 \operatorname{drift}\bigl(\dd F_t(s)\bigr)
 &=\Bigl[
 \frac12s^2\partial_{ss}F_t(s)
 +\frac{\theta}{2}s\partial_sF_t(s)
 +\frac{\kappa}{2}sG_t(s)f'(G_t(s))\\
 &-\frac12s^2
 \bigl(2f'(G_t(s))+G_t(s)f''(G_t(s))\bigr)
 \partial_{ss}G_t(s)
 \Bigr]\,\dd t.
\end{aligned}
\]
Here the derivatives of $f$ are with respect to its scalar argument.
For a general $f$, taking expectations (assuming the required integrability is valid) introduces additional
quantities through the last two terms.
For example, $f(u)=u^{-2}$ leaves the remainder
$-s^2G_t(s)^{-3}\partial_{ss}G_t(s)$, whose expectation must also
be determined to describe the evolution of $\E[F_t(s)]$.

To obtain closure with the spatial operator
$\mathcal G_1^\theta=\frac12s^2\partial_{ss}
+\frac{\theta}{2}s\partial_s$,
the derivative remainder vanishes when
\[
 uf''(u)+2f'(u)=0.
\]
The solutions are $f(u)=c_0+c_1/u$.
Choosing the reciprocal part $f(u)=u^{-1}$ also turns the term
involving $\kappa$ into multiplication by $-\kappa s/2$.
Thus, if we choose $h_t(s)=G_t(s)^{-1}$, we have
\[
 \operatorname{drift}\bigl(\dd h_t(s)\bigr)
 =
 \left(\mathcal G_1^\theta-\frac{\kappa}{2}s\right)
 h_t(s)\,\dd t.
\]
The resulting operator acts only on the probe argument $s$, with
coefficients independent of the state $\boldsymbol{\Phi}_t$.
The bounds and martingale argument in the preceding proof therefore
give a closed linear equation for the corresponding event-weighted
expectations.

The diffusion interpretation provides a Feynman--Kac representation
of this equation.  In particular, for a deterministic state $g$, put
\[
 \Psi_{\tau,s}(g)
 =
 \E_s\left[
 \exp\left(-\frac{\kappa}{2}\int_0^\tau
                     \mathsf S_v\,\dd v\right)
 \frac1{g(-\mathsf S_\tau)}
 \right],
\]
where $\mathsf S$ has generator $\mathcal G_1^\theta$.
The conditional duality identity is then (recall the notation $h_g^+$)
\[
 \E\left[
 h_{\boldsymbol{\Phi}_t}^+(s)\mid\mathscr F_r
 \right]
 =
 \Psi_{t-r,s}(\boldsymbol{\Phi}_r).
\]
The dual law and its weight are fixed independently of the original
solution.  Thus \(\Psi_{\tau,s}\) gives a solution-independent
prescription for the conditional expectation of each one-point
probe, depending on the past only through
\(\boldsymbol{\Phi}_r\).

To determine the full conditional law, we use all finite products
of reciprocal probes.  In It\^o's product rule, the
cross-covariations become the pairwise first-order interaction
terms in the finite-dimensional dual generators, so the product
observables also satisfy closed equations.
Analogous arguments apply in the Hua--Pickrell case.
For each model, taking all finite products at rational probe
arguments, together with the constant function \(1\), gives
the bounded determining class used in
Proposition~\ref{prop:conditional-determining-algebra}.
The conditional dualities for these classes therefore yield
uniqueness in law.
\end{rmk}

\subsection{Determining classes and uniqueness}

Conditional identities on a countable determining class are enough to
construct a Borel transition kernel and identify the entire continuous path law.
\begin{prop}
\label{prop:conditional-determining-algebra}
Let ${\mathscr E}$ be Polish and let
$\mathscr M\subset C_b({\mathscr E};\C)$ be countable, unital,
multiplicatively closed, self-adjoint, and separating.
For each $t\geq0$ and $F\in\mathscr M$, fix a bounded Borel function
$\mathcal P_tF$ on ${\mathscr E}$.
Suppose that, for every $x\in{\mathscr E}$, there exists a continuous
adapted ${\mathscr E}$-valued process $\mathsf X$ on a filtered
probability space, with $\mathsf X_0=x$ almost surely, satisfying
\begin{equation}
 \label{eq:abstract-conditional-duality}
 \E\left[F\left(\mathsf X_t\right)\mid\mathscr F_r\right]
 =\left(\mathcal P_{t-r}F\right)\left(\mathsf X_r\right),
 \qquad F\in\mathscr M,\quad 0\leq r\leq t.
\end{equation}
Then there exists a Borel Markov semigroup such that every continuous
adapted ${\mathscr E}$-valued process satisfying
\eqref{eq:abstract-conditional-duality} is Markov with this semigroup.
Any two such processes with the same initial law have the same law
on $C(\R_+;{\mathscr E})$.
\end{prop}

\begin{proof}
Enumerate $\mathscr M=(F_j)_{j\geq1}$.  The continuous injection
$x\mapsto(F_j(x))_j$ into $\C^\N$ and Lusin--Souslin give
$\sigma(\mathscr M)=\mathscr B({\mathscr E})$.  The map
\[
 \iota:\mathscr P\left({\mathscr E}\right)\longrightarrow\C^\N,\qquad
 \iota\left(\rho\right)=\left(\int_{\mathscr E}F_j\,\dd\rho\right)_j
\]
is Borel and injective: equality of its coordinates gives equality on the
unital self-adjoint algebra spanned by $\mathscr M$, and then on bounded
Borel functions by the functional monotone-class theorem.  Lusin--Souslin
makes its range and inverse Borel.  Existence and
\eqref{eq:abstract-conditional-duality} at $r=0$ show that
$\mathbf v_t(x)=((\mathcal P_tF_j)(x))_j$ lies in this range, so
$\mathbf P_t(x,\cdot)=\iota^{-1}(\mathbf v_t(x))$ is a Borel probability kernel.
The conditional identity extends to every bounded Borel function by
linear combinations and bounded pointwise convergence, using conditional
dominated convergence and the functional monotone-class theorem.
Successive conditioning gives the semigroup law and all finite-dimensional
distributions.  Rational evaluations generate the Borel sigma-field of
$C(\R_+;{\mathscr E})$, proving path-law uniqueness.
\end{proof}

\begin{proof}[Proof of Theorem~\ref{thm:intro-uniqueness}]
Existence follows from
\cite[Theorems~1.4 and~1.6]{AssiotisMirsajjadi}.
We verify the hypotheses of
Proposition~\ref{prop:conditional-determining-algebra}.
Set
$\mathscr Q_+=\mathbb Q_{>0}$,
$\mathscr Q_{\mathrm{HP}}=\mathbb Q\setminus\{0\}$,
and, for $\diamond\in\{+,\mathrm{HP}\}$, let
\[
 \mathscr M_\diamond\defeq\{1\}\cup
 \left\{
 H_n^\diamond(\cdot,\mathbf q):
 n\geq1,\ \mathbf q\in\mathscr Q_\diamond^n
 \right\}.
\]

First suppose that $\mathbf q$ has pairwise distinct coordinates and
write $F=H_n^\diamond(\cdot,\mathbf q)$.
Define $(\mathcal P_tF)(g)$ by the corresponding dual expectation in
\eqref{eq:conditional-duality} or
\eqref{hp:eq:conditional-duality}, with $r=0$ and
$\boldsymbol{\Phi}_0$ replaced by the deterministic state $g$.
Thus the expectation defining $\mathcal P_tF$ is over the dual alone.
Unordered tuples are treated by increasing rearrangement, since the
probes and the dual weights are invariant under relabelling.
Put $\mathcal P_t1=1$.
The probes and the exponential weights have modulus at most one, so
$|\mathcal P_tF|\leq1$; dominated convergence under the fixed dual law
also shows that $\mathcal P_tF$ is continuous in $g$.

Proposition~\ref{prop:conditional-duality} now gives, for every solution
$\boldsymbol{\Phi}$ with any initial law,
\[
 \E\left[
 F(\boldsymbol{\Phi}_t)\mid\mathscr F_r
 \right]
 =
 (\mathcal P_{t-r}F)(\boldsymbol{\Phi}_r),
 \qquad 0\leq r\leq t,
\]
whenever $F$ is a probe with distinct arguments.
This is precisely the conditional identity required in
\eqref{eq:abstract-conditional-duality}, initially for this subclass.

We next extend the definition and the identity to repeated arguments.
Given $\mathbf q\in\mathscr Q_\diamond^n$, choose
$\mathbf q^{(m)}\to\mathbf q$ with pairwise distinct coordinates and put
\[
 F_m=H_n^\diamond(\cdot,\mathbf q^{(m)}),
 \qquad F=H_n^\diamond(\cdot,\mathbf q).
\]
Then $F_m(g)\to F(g)$ for every state $g$, and $|F_m|\leq1$.
Fix $g$ and choose any solution started from $g$.
Applying Proposition~\ref{prop:conditional-duality} at $r=0$ and taking
expectations gives
\[
 (\mathcal P_tF_m)(g)
 =
 \E_g\left[F_m(\boldsymbol{\Phi}_t)\right]
 \longrightarrow
 \E_g\left[F(\boldsymbol{\Phi}_t)\right],
\]
where $\E_g$ denotes expectation for the chosen solution and the
convergence follows by bounded convergence.
Consequently, we may define
\[
 (\mathcal P_tF)(g)
 \defeq\lim_{m\to\infty}(\mathcal P_tF_m)(g).
\]
The preceding argument shows that this limit exists and is independent
of the approximation, and that every solution from $g$ gives the same
value.  Only existence of a solution from $g$ is used here.
The same argument identifies different tuples representing the same
function.  Hence $\mathcal P_tF$ is unambiguous, bounded by one in
modulus, and Borel as a pointwise limit of Borel functions.

For an arbitrary solution and $0\leq r\leq t$, bounded convergence gives
\[
 F_m(\boldsymbol{\Phi}_t)\longrightarrow F(\boldsymbol{\Phi}_t),
 \qquad
 (\mathcal P_{t-r}F_m)(\boldsymbol{\Phi}_r)
 \longrightarrow
 (\mathcal P_{t-r}F)(\boldsymbol{\Phi}_r)
 \quad\text{in }\mathrm L^1.
\]
The $\mathrm L^1$ contraction of conditional expectation therefore
allows us to pass to the limit in the conditional identities supplied
by Proposition~\ref{prop:conditional-duality}.  We obtain
\[
 \E\left[
 F(\boldsymbol{\Phi}_t)\mid\mathscr F_r
 \right]
 =
 (\mathcal P_{t-r}F)(\boldsymbol{\Phi}_r),
 \qquad F\in\mathscr M_\diamond,\quad 0\leq r\leq t.
\]

Each class $\mathscr M_\diamond$ consists of bounded continuous
functions and is countable and unital.
It is multiplicatively closed because multiplication concatenates
the argument tuples, with repetitions allowed.
It is self-adjoint because positive probes are real and
$h_g^{\mathrm{HP}}(-q)=\overline{h_g^{\mathrm{HP}}(q)}$.
If two states have equal one-point probes, their entire functions
agree on the negative rational axis or the nonzero imaginary rational
axis; the identity theorem therefore gives separation.
Consequently,
\begin{equation}
 \label{eq:algebra-generates}
 \sigma(\mathscr M_+)=\mathscr B(\mathscr E_+),
 \qquad
 \sigma(\mathscr M_{\mathrm{HP}})=\mathscr B(\mathscr E).
\end{equation}
Since the state spaces are Polish and solutions exist from every
deterministic state, all the hypotheses of
Proposition~\ref{prop:conditional-determining-algebra} are satisfied.
That proposition proves uniqueness in law, in fact for every initial
law, and supplies the Borel Markov semigroups.
\end{proof}
\subsection{Relation to the boundary intertwinings}
\label{sec:boundary-intertwinings}

The reciprocal dualities are related directly to the intertwinings used
to construct the boundary processes in
\cite{AssiotisHP,AssiotisMirsajjadiISDE,AssiotisMirsajjadi} in the first place. We now explain how.

Identify the boundaries with the entire-function spaces of
Section~\ref{sec:state-spaces}.  For either model, write
\(\mathbf P_t^\diamond\) for its boundary semigroup,
\(\mathbf P_{n,t}^\diamond\) for the original finite-dimensional semigroup,
and \(\Lambda_n^\infty\) for the boundary-to-corner eigenvalue kernel \cite{AssiotisHP,AssiotisMirsajjadiISDE}.
Let \(\overline{\mathscr W}_n^\diamond\) denote the corresponding closed
positive or real Weyl chamber, with increasing coordinates.
The boundary construction gives
\(\mathbf P_t^\diamond\Lambda_n^\infty
=\Lambda_n^\infty\mathbf P_{n,t}^\diamond\), for every \(n\) and \(t\geq0\).

Let \(\mathbf D_{\mathbf v}=\operatorname{diag}(v_1,\ldots,v_n)\), put
\(\varepsilon_+=-1\), \(\varepsilon_{\mathrm{HP}}=\mathrm i\), and define
the orbital Laplace and Fourier kernels, with normalized Haar measure, by
\[
 \mathcal K_n^\diamond\left(\mathbf x,\mathbf v\right)
 =\int_{\mathscr U\left(n\right)}
 \mathrm e^{\varepsilon_\diamond
   \Tr\left(\mathbf D_{\mathbf v}\mathbf U\mathbf D_{\mathbf x}\mathbf U^*\right)}\,\dd\mathbf U.
\]
Here \(\mathbf v=\mathbf s\in\mathscr W_n\) in the positive case and
\(\mathbf v=\mathbf u\) belongs to a signed chamber in the Hua--Pickrell
case.  For fixed \(g\), the Gaussian matrix representation of \(\Lambda_n^\infty\)
\cite[Definitions~2.17--2.18]{AssiotisMirsajjadiISDE} has independent
diagonal entries in its corner, with common law that of
\(\mathsf Y_g^+\) in \eqref{eq:Laplace-representation}, or
\(\mathsf Y_g\) in Lemma~\ref{hp:lem:characteristic-representation}.
Their transforms and unitary invariance therefore give
\[
 H_n^\diamond\left(g,\mathbf v\right)
 =\int_{\overline{\mathscr W}_n^\diamond} \mathcal K_n^\diamond\left(\mathbf x,\mathbf v\right)
       \Lambda_n^\infty\left(g,\dd\mathbf x\right).
\]
Thus the reciprocal product is the joint transform of the diagonal
entries of the corner.  The orbital average is essential when working
with its eigenvalues.  Both kernels have modulus at most one, so the
identity holds under arbitrary boundary laws without moment assumptions.

Write \(\mathcal D_n^\diamond\) for the generator of
\(\mathbf P_{n,t}^\diamond\).  Conjugating it by the Vandermonde reduces
it, up to a constant, to a sum of the scalar operators
\[
 \begin{aligned}
 &\frac{x^2}{2}\partial_x^2+
   \left[\left(\frac\theta2-n+1\right)x+\frac\kappa2\right]\partial_x,
 &&\diamond=+,\\
 &\left(1+x^2\right)\partial_x^2+
   \left[\left(2-2n-2a\right)x+2b\right]\partial_x,
 &&\diamond=\mathrm{HP}.
 \end{aligned}
\]
The HCIZ formula expresses \(\mathcal K_n^\diamond\) as a constant times
\(\det[\mathrm e^{\varepsilon_\diamond x_iv_j}]/
  [\Delta(\mathbf x)\Delta(\mathbf v)]\).
Expanding the determinant, the quadratic and linear homogeneous terms
exchange \(\mathbf x\) and \(\mathbf v\).  The remaining terms give
\(-\kappa\sum_{i=1}^{n} s_i/2\) and
\(2\mathrm i b\sum_{i=1}^{n} u_i-\sum_{i=1}^{n} u_i^2\), respectively.  Undoing the
Vandermonde conjugation cancels the dimension shifts and gives
\[
 \mathcal D_{n,\mathbf x}^\diamond \mathcal K_n^\diamond
 =\left(\mathcal B_{n,\mathbf v}^\diamond+C_\diamond\left(\mathbf v\right)\right)
    \mathcal K_n^\diamond,\qquad
 \begin{cases}
 \mathcal B_n^+=\mathcal G_n^\theta,&C_+=-\kappa V_n^+,\\
 \mathcal B_n^{\mathrm{HP}}=\mathcal A_n^a,&C_{\mathrm{HP}}=C_{b,n}.
 \end{cases}
\]

This generator identity lifts to the finite semigroups.  From deterministic
finite initial states the positive trace and the Hua--Pickrell squared
norm have all polynomial moments on bounded time intervals.
Derivatives of the orbital kernels have polynomial growth, and symmetry
removes the apparent collision poles.  Localized It\^o formulas and these
moment bounds justify the finite duality.  Feller continuity extends it to closed-chamber
initial states.  Writing \(\mathcal T_{n,t}^\diamond\) for the dual
Feynman--Kac operator, this gives
\[
 \mathbf P_{n,t}^\diamond \mathcal K_n^\diamond\left(\cdot,\mathbf v\right)
       \left(\mathbf x\right)
 =\mathcal T_{n,t}^\diamond \mathcal K_n^\diamond\left(\mathbf x,\cdot\right)
       \left(\mathbf v\right).
\]
Consequently, composing with the boundary intertwining yields the duality
\[
 \mathbf P_t^\diamond H_n^\diamond\left(\cdot,\mathbf v\right)\left(g\right)
 =\Lambda_n^\infty\mathbf P_{n,t}^\diamond
       \mathcal K_n^\diamond\left(\cdot,\mathbf v\right)\left(g\right)
 =\mathcal T_{n,t}^\diamond H_n^\diamond\left(g,\cdot\right)\left(\mathbf v\right).
\]
Conversely, the orbital transforms determine finite eigenvalue laws:
their unitary-invariant matrix lifts are determined by their Laplace or
Fourier transforms.  Hence the reciprocal dualities, for all \(n\),
recover the boundary intertwinings for any specified boundary Markov
semigroup.

\bibliographystyle{amsplain}
\bibliography{ReferencesLP}
\bigskip

\noindent{\sc School of Mathematics, University of Edinburgh, James Clerk Maxwell Building, Peter Guthrie Tait Rd, Edinburgh EH9 3FD, U.K.}\newline
\href{mailto:theo.assiotis@ed.ac.uk}{\small theo.assiotis@ed.ac.uk}

\bigskip

\noindent{\sc School of Mathematics, University of Edinburgh, James Clerk Maxwell Building, Peter Guthrie Tait Rd, Edinburgh EH9 3FD, U.K.}\newline
\href{mailto:F.Li-31@sms.ed.ac.uk}{\small F.Li-31@sms.ed.ac.uk}

\end{document}